\documentclass{article}

\def\n{\noindent}
\def\Z{{\bf Z}}
\def\N{{\bf N}}
\def\Q{{\bf Q}}
\def\R{{\bf R}}
\def\F{{\bf F}}
\def\C{{\bf C}}

\begin{document}

\hspace{2.3 in} {\it to appear in Izvestiya Math.} {\bf 77} (2013)
  
\bigskip

\bigskip

\bigskip

\centerline
{\bf Hasse principle for $G$--trace forms}

\bigskip

\bigskip
\qquad \qquad \qquad {E. Bayer--Fluckiger, R. Parimala, J--P. Serre}

\bigskip

\bigskip

\centerline {\it \`A Igor Rostislavovich Shafarevich, \`a l'occasion de son $90$-i\`eme anniversaire}

\medskip
\medskip
{\bf Introduction}

\bigskip
 
 Let  $k$  be a field, let $G$ be a finite group and let $L$ be a $G$-Galois algebra over $k$ , for instance a finite Galois extension of  $k$ with Galois group $G$. The {\it trace form} of $L$ is the symmetric bilinear form
 
\smallskip
\quad \quad \quad $q_L : L \times L \to k$,  \  \ where \ \ $q_L(x,y) = {\rm Tr}_{L/k}(xy).$

\smallskip 

Note that $(L,q_L)$ is a {\it $G$--quadratic space}, 
i.e. we have $q_L(gx,gy) =
q_L(x,y)$ for all 
$g \in G$ and all $x,y \in L$. 

\smallskip

An interesting special case is when $L$ has a {\it self-dual normal basis},
i.e. when there exists $x \in L$ such that $(gx)_{g\in G}$ is a $k$-basis of $L$ and $q_L(gx,hx) = \delta_{g,h}$ for every $g,h \in G$, where $\delta_{g,h}$ is the Kronecker symbol. 

\smallskip

The aim of the present paper is to prove (cf. th.1.3.1):

\smallskip
\n
{\bf Theorem} - {\it Suppose that $k$ is a global field of characteristic $\neq 2$. Let $L$ and $L'$
be two $G$--Galois algebras. If the $G$--quadratic spaces  $(L,q_L)$ and $(L',q_{L'})$ become isomorphic over
all the completions of $k$, then they are  isomorphic over~$k$.}

\smallskip
\n
{\bf Corollary} - {\it If $L$ has a self--dual normal basis over all the completions of $k$, then it has a self-dual normal basis over $k$.}

\n [In short: the Hasse principle applies to self-dual normal bases in characteristic $\neq 2$.]

\smallskip

These results complete those of  [BFP 11] 
where it was assumed that  $G$ has the ``odd
determinant property''. The present proof follows a similar pattern. We give it in \S1; the only new ingredient is a restriction-induction property which is proved in \S4 by using Burnside rings, following a method introduced by A. Dress in the 1970s. There are also two  basic results which have to be used:

\smallskip

\n $\bullet$ Witt cancellation theorem for $G$-quadratic spaces; 
 
 \smallskip
\n $\bullet$ the ``odd does not count" principle, which allows division by any odd integer.

\smallskip
   Their proofs are summarized  (resp. completed) in \S2 (resp. in \S3), in the more general setting of hermitian spaces over an algebra with involution in characteristic $\neq 2$.

     \medskip

\n {\it Acknowledgements}. We thank  H.W. Lenstra, Jr, and J-P. Tignol for their help with \S3.4 and  \S3.5.

\medskip

\smallskip

\medskip
\medskip
{\bf \S 1. $G$--trace forms and Hasse principle}

\medskip
The aim of the present \S \ is to prove the Hasse principle result stated in the introduction, assuming the technical results of \S\S2,3,4.

Before doing so, we introduce some notation and basic facts concerning $G$--quadratic
forms and $G$--trace forms.

\smallskip
 
In what follows, the characteristic of the ground field $k$ is assumed to be $\neq 2$.

\medskip

\medskip
{\bf 1.1. $G$--quadratic spaces}

\medskip
Let  $G$ be a finite group. 
A {\it $G$--quadratic space} is a pair $(V,q)$, where $V$ is a $k[G]$--module
that is a finite dimensional $k$--vector space, 
and $q : V \times V \to k$ is a non--degenerate symmetric bilinear
form such that $$q(gx,gy) = q(x,y)$$ for all $x,y \in V$ and all $g \in G$. 

Two $G$--quadratic spaces  $(V,q)$ and $(V',q')$ are {\it isomorphic} if there exists
an isomorphism of $k[G]$--modules $f : V \to V'$ such that $q'(f(x),f(y)) = q(x,y)$ for
all $x,y \in V$. If this is the case, we write $(V,q) \simeq_G (V',q')$.

 \medskip
\n
{\bf Theorem 1.1.1}  ([BFL 90, th. 4.1]) - {\it If  two $G$--quadratic spaces
become isomorphic over an odd degree extension, then they are isomorphic.}

\medskip
Let $S$ be a subgroup of $G$. We have two operations, induction and restriction (see
for instance [BFS 94, \S1.2]) :

\medskip
\n $\bullet$ If $(V,q)$ is an $S$--quadratic form, then ${\rm Ind}_S^G(V,q)$ is a $G$--quadratic form;

\smallskip
\n $\bullet$ If $(V,q)$ is a $G$--quadratic form, then ${\rm Res}_S^G(V,q)$ is an $S$--quadratic form.

\medskip
The following result is one of the basic tools in the proof of the Hasse principle~:

\medskip
\n
{\bf Theorem 1.1.2} - {\it Let $S$ be a $2$--Sylow subgroup of $G$, and let $(V_1,q_1)$ and $(V_2,q_2)$ be two $S$--quadratic spaces over $k$. Suppose that 

$${\rm Res}^G_S \  {\rm Ind}^G_S (V_1,q_1) \simeq_S {\rm Res}^G_S \ {\rm Ind}^G_S (V_2,q_2).$$

\n Then

$${\rm Ind}^G_S (V_1,q_1) \simeq_G {\rm Ind}^G_S (V_2,q_2).$$}

\medskip
 This was proved in [BFS 94, th.5.3.1] when $S$ is elementary abelian and  in [BFP 11, theorem 2.2] when the fusion of $S$ in $G$ is controlled by the normalizer of $S$. The proof in the general case will be given in \S 4.

\medskip
Note that the Hasse principle does not hold in general for $G$--quadratic spaces (see e.g.  [M 86, \S 2]). However,
it does hold when $G$ is a 2--group :

\medskip
 \n
 {\bf Proposition 1.1.3} - {\it Suppose that $k$ is a global field, and let $S$ be a $2$--group. If two $S$--quadratic spaces become isomorphic over all
 the completions of $k$, then they are isomorphic over $k$}.
  
 \medskip
 
 \n {\it Proof.} By a theorem of M. Kneser (see e.g. [R 11, th.3.3.3]), this is true for any finite group  $S$, of order prime to the characteristic of $k$, which has the following property:
 
 \medskip
 
 (P) - $k[S]$ does not have any simple factor, stable under the involution, which is a matrix algebra over a quaternion skewfield, the involution being of orthogonal type. 
 
 \medskip
 
\n  By [BFP 11, prop.3.7], property (P) holds when $S$ is a 2-group and the characteristic of $k$ is not 2.

\bigskip

{\bf 1.2. Trace forms}

\medskip
Let $L$ be a $G$--Galois algebra over $k$. Recall (cf. e.g. [BFS 94, \S1.3]
or [A8, \S16.7]) that this means that $L$ is an \'etale $k$-algebra on which $G$ acts by $k$-automorphisms in such a way that $L$ is a free $k[G]$-module of rank 1. Let 
$$q_L : L \times L \to k, \ \ q_L(x,y) = {\rm Tr}_{L/k}(xy)$$
be the trace form of $L$. Then $(L,q_L)$ is a $G$--quadratic space.  

\medskip
\n
{\bf Lemma 1.2.1} - {\it Let $S$ be a $2$--Sylow subgroup of $G$. If $L$ is a $G$--Galois
algebra, there exists an odd degree field extension $k'/k$ and an $S$--Galois algebra
$M$ over $k'$ such that the $G$--form $(L,q_L ) \otimes_k k'$ is isomorphic to the
$G$--form ${\rm Ind}_S^G(M,q_M)$.}

\smallskip

\n This is proved in [BFS 94, prop.2.1.1].

\medskip
{\bf 1.3. Proof of the Hasse principle}

\medskip
We now prove the result stated in the introduction:

\medskip
\n
{\bf Theorem 1.3.1} - {\it Suppose that $k$ is a global field. Let $L_1$ and $L_2$
be two $G$--Galois algebras. If the $G$--quadratic spaces  $(L_1,q_{L_1} \! )$ and $(L_2,q_{L_2} \! )$ become isomorphic over
all the completions of $k$, then they are  isomorphic over $k$.}

\medskip
\n
{\it Proof.} 
By lemma 1.2.1, there exists an odd degree field extension $k'/k$ and two $S$--Galois algebras
$M_1$ and $M_2$ over $k'$ such that 
$$(L_1,q_{L_1} ) \otimes_k k' \simeq_G {\rm Ind}_S^G(M_1,q_{M_1}) \ \ {\rm and} \  \ (L_2,q_{L_2})\otimes_k k' \simeq_G {\rm Ind}_S^G(M_2,q_{M_2}).$$  

Suppose that the $G$--quadratic spaces 
$(L_1,q_{L_1} ) $ and $ (L_2,q_{L_2})$ are isomorphic over all the completions of $k$.  This implies
that the $G$--quadratic spaces  $$(L_1,q_{L_1} )\otimes_k k' \quad {\rm and} \quad 
 (L_2,q_{L_2}) \otimes_k k' $$ are isomorphic over all the completions of $k'$.  Therefore the 
$S$--quadratic spaces
$${\rm Res}^G_S (L_1,q_{L_1} ) \otimes_k k'  \simeq_S {\rm Res}^G_S\ {\rm  Ind}^G_S (M_1,q_{M_1})$$ and $$ {\rm Res}^G_S  (L_2,q_{L_2}) 
\otimes_k k' \simeq_S {\rm Res}^G_S\ {\rm  Ind}^G_S(M_2,q_{M_2})$$  are isomorphic over all the completions of $k'$.
 Since the Hasse principle holds for $S$--quadratic spaces (cf. prop.1.1.3), 
this implies that $${\rm Res}^G_S\ {\rm  Ind}^G_S (M_1,q_{M_1}) 
\simeq_S {\rm Res}^G_S\ {\rm  Ind}^G_S(M_2,q_{M_2}).$$

\n By th.1.1.1 we have

$${\rm  Ind}^G_S (M_1,q_{M_1}) 
\simeq_G {\rm  Ind}^G_S(M_2,q_{M_2}).$$ Recall that we have
$$(L_1,q_{L_1} ) \otimes_k k' \simeq_G {\rm Ind}_S^G(M_1,q_{M_1}) \quad {\rm and} \quad  (L_2,q_{L_2})\otimes_k k' \simeq_G {\rm Ind}_S^G(M_2,q_{M_2}),$$  hence 
$$(L_1,q_{L_1} ) \otimes_k k' \simeq_G (L_2,q_{L_2})\otimes_k k'.$$ By th.1.1.1, we
get $$(L_1,q_{L_1} ) \simeq_G (L_2,q_{L_2}),$$ as claimed.

\medskip
\medskip

\centerline
{\bf  \S 2. Witt's cancellation theorem for $\epsilon$--hermitian spaces}

\medskip
 Witt's cancellation theorem applies not only to quadratic forms, but also
to hermitian spaces over finite dimensional algebras with
involution, provided the characteristic of the ground field is $\not = 2$. This is due to H-G. Quebbemann, R. Scharlau, W. Scharlau and M. Schulte, cf.  [QSSS 76, 3.4 (i)]; see also [QSS 79, 3.4.(iii)] and [K 91 chap.II, th.6.6.1]. In what follows we fix the notation and we summarize the different steps of the proof.

\medskip 

{\bf 2.1. $\epsilon$--hermitian spaces }

\smallskip From now on, unless otherwise stated, all modules are left modules.

\smallskip

Let $R$ be a ring endowed with an involution $r \mapsto \overline r$. For any 
$R$--module $M$, we denote by $M^*$ its dual  ${\rm Hom}_R(M,R)$. Then $M^*$  has an $R$--module structure given by $(rf)(x) =  f(x) \overline r$ for all $r \in R$, $x \in M$ and
$f \in M^*$. If $M$ and $N$ are two $R$--modules and if $f : M \to N$ is
a homomorphism of $R$--modules, then $f$ induces a homomorphism
$f^* : N^* \to M^*$ defined by $f^*(g) = g \circ f$ for all $g \in N^*$. 

\smallskip
Let $\epsilon = \pm 1$ and let $M$ be an $R$-module. An $\epsilon$--{\it hermitian form on $M$} is a biadditive map $h : M \times M \to R$ which satisfies the following condition:
\smallskip
  $h(rx,sy) = r h(x,y) \overline s$ \ and \ $\overline {h(x,y)} = \epsilon h(y,x)$ \ for
all $x,y \in M$ \ and all $r,s \in R$.

\smallskip

\n In that case, the map $ h : M \to M^*$ given by  $y \mapsto (x \mapsto (h(x,y))$ is an $R$--homomorphism. If it is an isomorphism, we say that $h$ is {\it nonsingular}, and we call the pair $(M,h)$ an {\it $\epsilon$--hermitian space}. 

\smallskip

  If $(M,h)$ and $(M',h')$ are two $\epsilon$--hermitian spaces, their orthogonal sum is also an $\epsilon$--hermitian space; it is denoted by  $(M,h) \ \oplus \ (M',h')$.
  
  \smallskip
  
  \n {\it Remark on notation}
  
    If $(M,h)$ is an  $\epsilon$--hermitian space, we are using the letter $h$  for two different, but closely related, objects: a biadditive map
 $M \times M \to R$, and an $R$--linear homomorphism $M \to M^*$.
  There is usually no risk of confusion, but, when there is one, we shall write
  $h_{{\rm b}}$ for the biadditive map \ $M \times M \to R$ \ and keep $h$ for the map $M \to M^*$.
  
  \smallskip
  
  \n {\it Example $:$ the group algebra $k[G]$}
  
  Let $k$ be a commutative ring, and let $G$ be a finite group. The group
  algebra  $R = k[G]$ has a natural $k$--linear involution, characterized by
  the formula $\overline g = g^{-1}$ for every $g \in G$. We have the following dictionary:
  
   a) $R$--module $M  \ \Longleftrightarrow  k$-module $M$ with a $k$--linear action of $G$;
   
   b) $R$--dual $M^*\ \Longleftrightarrow k$--dual of $M$, with the contragredient (i.e. dual) action of $G$. 
   
  c) $\epsilon$--hermitian space $(M,h)$\ $\Longleftrightarrow$ 
  $\epsilon$--symmetric bilinear form on $M$, which is $G$--invariant and defines an isomorphism of $M$ onto its $k$--dual.
  
  \smallskip
  
  \n [Let us be more specific about these identifications. Denote by $\alpha: k[G] \to k$ the linear form ``coefficient of  1'', i.e. $\alpha(g)=1$ if $g=1$ and $\alpha(g) = 0$ if $g\neq 1$. The identification of b) transforms $y\in M^* = {\rm Hom}_R(M,R)$ into the $k$--linear form  $\alpha \circ y$. Similarly, the identification of  c)  transforms $h_{{\rm b}}$ into the $k$--bilinear form $(x,y) \mapsto \alpha(h_{\rm b}(x,y))$.] 
  
  \smallskip
  
    Note that, when $k$ is a field, condition c) implies that $M$ is finite dimensional over $k$ (if not it would not be isomorphic to its dual), so that,
    when $\epsilon = 1$, an $\epsilon$--hermitian space over $k[G]$ is the same as what we called a $G$--quadratic space in \S1.1.

\smallskip

\medskip

{\bf 2.2. Basic hypotheses for \S2 and \S3}

\smallskip

 In the rest of \S2 and in  \S3,  we assume that $R$ is an algebra over a field $k$ and that its involution $r \mapsto \overline r$ is $k$--linear. We also assume:
 
   (2.2.1)  \ {\it The characteristic of $k$ is $\neq 2$}.
   
   (2.2.2) \  {\it ${\rm dim}_k \ R < \infty.$}

   \smallskip
   
   \n[The only case we actually need for theorem 1.3.1 is  $R = k[G]$, with $G$ finite.]  
    \smallskip
   
   Note that  (2.2.2)  implies that $R$ is a left and right artinian ring, cf. [A8, \S1]. 
   
   \smallskip
    
\n    If $(M,h)$ is an $\epsilon$--hermitian space over $R$, we shall always assume:
    
    \smallskip
    
    (2.2.3) \ {\it ${\rm dim}_k M < \infty$.}

    \smallskip
    
    Because of (2.2.2), such an $M$ is an $R$--module of finite length, so that one may apply the Krull-Remak-Schmidt theory to it;  in particular the cancellation theorem holds, cf. [A8, \S2.4, cor.3 to th.1]. 
   
       \medskip

{\bf 2.3.  Biduality and reflexivity}

\smallskip

  If $M$ is an $R$-module, the dual $(M^*)^*$ of its dual is called its {\it bidual}
  and is denoted by $M^{**}$. Every element  $x\in M$ defines an $R$-linear form $\ell_x$ on $M^*$ by the formula  $\ell_x(f) = \overline{f(x)}$, for $f \in M^*$. The map $x \mapsto \ell_x$ is an $R$-linear homomorphism  $  M \ \to \ M^{**},$   which we shall denote by $c_M$. The same construction applies to $M^*$, so that we have $c_{M^*}: M^* \to M^{***}$. On the other hand, the adjoint $(c_M)^*$ of $c_M$ is a homomorphism of $M^{***}$ into $M^*$.
 These two maps are related by :
 
 \smallskip
 
 \n {\bf Proposition 2.3.1} - {\it The composition
 
 $$ (c_M)^* \circ  c_{M^*} \ : \ M^* \  \to \  M^{***} \ \to \ M^*$$
 
\n is the identity.}

\smallskip

The proof is a direct computation. To explain it, it is convenient to
denote the duality pairing between  $M$ and  $M^*$ by a scalar product,
i.e. to write
 $ <\!x,y\!>$ instead of $y(x)$, for  $x\in M$ and $y \in M^*$. With a similar convention for  the other modules, one finds:
 
 \smallskip
 
 $ <\!x, (c_M)^* \circ  c_{M^*}(y)\!> \ = \ <\!c_M(x),c_{M^*}(y)\!> \ = \ <\!c_M(x),y\!> \ = \ <\!x,y\!> ,$
 
 \smallskip
\n and since this is valid for every $x$, we have $(c_M)^* \circ  c_{M^*}(y)=y$ for every $y \in M^*$, as claimed.

\smallskip
 
 A module $M$ is said to be {\it reflexive} if  $c_M :   M \ \to M^{**}$ is an isomorphism; it is said to be {\it selfdual} if $M$ and $M^*$ are isomorphic.
 
 \smallskip
 
\n {\bf Proposition 2.3.2} - (1) {\it If $M \simeq M_1 \oplus M_2$, then $M$
is reflexive if and only if both $M_1$ and $ M_2$ are reflexive.}

(2) {\it If a module $N$ is isomorphic to the dual of some module, then $c_N$ is injective and its image is a direct factor in $N^{**}.$}

(3) {\it If $M$ is selfdual, $M$ is reflexive.}

(4) {\it If $M$ is reflexive, so is $M^*$.}
\medskip

\n {\it Proof.} Assertion (1) is immediate. Assertion (2) follows from prop.2.3.1 applied to a module $M$ such that  $N \simeq M^*$. In case (3),
 $M^*$ is isomorphic to $M$;
hence $M^{**} \simeq M^* \simeq M$, and the modules $M$ and $M^{**}$ have the same length; since  $c_M: M  \to M^{**}$  is injective (by (2)), it is
an isomorphism. In case (4),  $c_{M^*} $ is injective because of (2), and since 
$M^* \simeq M^{**}$, the length argument above shows that
$c_{M^*} $ is bijective.

\medskip

\n {\bf Corollary 2.3.3} - {\it If $(M,h)$ is an $\epsilon$--hermitian space, 
then $M$ is reflexive, and so are all its direct factors.}

\smallskip

Since $M$ is selfdual, this follows from parts (3) and (1)  of prop.2.3.2.

\medskip 

{\bf 2.4.  Hyperbolic hermitian spaces}

\smallskip

 Let $\epsilon = \pm 1$, and let $N$ be a reflexive module. Define an 
 $\epsilon$--hermitian form $h_N$ on $N \oplus N^*$ by
 the biadditive map
 
 \smallskip
 $  ((x,y), (x',y'))  \ \  \mapsto \ \  y'(x) + \epsilon y(x') \quad {\rm for} \quad x,x' \in N \ \ {\rm and} \ \ y,y' \in N^*.$
 
 \smallskip
 \n Since $N$ is reflexive, $h_N$ is nonsingular. The hermitian space  $(N \oplus N^*,h_N)$  is denoted by $H_\epsilon(N).$ Note that  $H_\epsilon(N) \simeq H_\epsilon(N^*)$.
 
 \smallskip
  An $\epsilon$--hermitian space $(M,h)$ is called {\it hyperbolic} if there
  exists a reflexive module $N$ such that $(M,h) \simeq H_\epsilon(N)$. 
  
  \smallskip
  
  This is equivalent to asking that the module $M$ splits as a direct sum $M = N \oplus N'$ where  $N$ and $N'$ are {\it totally isotropic}, i.e. $h_{\rm b}(x,y)=0$ when $x,y \in N$ and when $x,y \in N'$.
  
  \smallskip
  
  \n {\it Example}. If $(M,h)$ is an $\epsilon$--hermitian space, then $(M,h)\oplus (M,-h)$ is hyperbolic: apply the criterion to $N= $ diagonal of $M \oplus M$ and  $N'= $ antidiagonal of $M \oplus M$.

    \smallskip
  
  \n {\bf Proposition 2.4.1} - {\it Let $(M,h)$ and $(M',h')$ be two $\epsilon$--hermitian hyperbolic modules. If $M \simeq M'$, then $(M,h) \simeq (M',h')$.}
  
  \smallskip
  
    \n {\bf Lemma 2.4.2} - {\it Let $N$ and $N'$ be two reflexive modules
    such that  $N \oplus N^* \simeq N' \oplus N'^*$. Then there exists two reflexive modules $P$ and $Q$ such that $N \simeq P \oplus Q$ and     $N' \simeq P \oplus Q^*$.}
    
    \smallskip
    
    \n {\it Proof of the lemma.} We use induction on the length $\ell $ of $N$. The case $\ell = 0$ is trivial.  If $\ell > 0$, we may decompose $N$ as $S \oplus N_0$, where $S$ is indecomposable and the length of $N_0$ is $< \ell$.
    By the Krull-Remak-Schmidt theorem, $S$ is isomorphic to a direct
    factor of either  $N'$ or $N'^*$. In the first case, we have  $N' \simeq  S \oplus N'_0$ for some reflexive module $N'_0$, hence
    
    \smallskip
    
     $S \oplus N_0  \oplus S^* \oplus  N_0^*  \simeq  S  \oplus  N'_0  \oplus S^*  \oplus (N'_0)^*$. 
     
     \smallskip
     
     \n By the cancellation theorem for modules, this implies  $ N_0   \oplus N_0^*    \simeq   N'_0   \oplus (N'_0)^*$, and we apply the induction assumption to the pair 
    $(N_0,N'_0)$. The second case follows from the first one, applied to the pair $(N,N'^*)$.
    
    \smallskip
    
    \n {\it Proof of proposition 2.4.1.} Suppose that $(M,h) \simeq H_\epsilon(N)$
    and $(M',h') \simeq H_\epsilon(N')$. We have  $N \oplus N^* \simeq N' \oplus N'^*.$ By the lemma above, we may write $N$ as $ P \oplus Q$ and $N' $ as $P \oplus Q^*$. We then have:
    
  \quad   $ (M,h) \simeq H_\epsilon(N)  \simeq H_\epsilon(P) \oplus H_\epsilon(Q) \simeq  H_\epsilon(P) \oplus H_\epsilon(Q^*) \simeq H_\epsilon(N') \simeq  (M',h'),$
    
    \n as claimed.
  \medskip

  {\bf 2.5.  The cancellation theorem}
  
  \smallskip
  
    Let ${\bf Ref}_R$ be the category of all reflexive finitely generated $R$--modules. This is an additive category, with a {\it duality functor} $M \mapsto M^*$, as in [K 91, chap.II, \S6.2].
    
    \smallskip
    
 \n   {\bf Proposition 2.5.1} - {\it The category ${\bf Ref}_R$ has the properties {\rm C.1},  {\rm C.2} and  {\rm C.3} of} [K 91, chap.II, \S5.2 and \S6.2].
 
 \smallskip
 
\n {\it Proof.} Property C.1 asks that every idempotent splits; this follows from prop.2.3.2 (1). 
 
 Property C.2 has two parts. First, every object $M$ should split  as a finite direct sum of indecomposable objects; this is true because $R$ is an artinian ring. Second, if  $M$ is indecomposable, the ring 
$A = {\rm End}(M)$ should be a (non necessarily commutative) local ring, i.e. $A$ should be $ \neq 0$ and, for all $a \in A$, either $a$ or $1-a$ should be invertible; this is true, cf. [A8, \S2.5, th.2].

Property C.3 asks that, for every object $M$, the ring $A = {\rm End}(M)$
should be complete for the topology defined by the powers of its radical;
indeed,  $A$ is artinian (because it is a finitely generated $k$--module), hence its radical is nilpotent ([A8, \S10.1, prop.1]), and the corresponding topology is the discrete topology.
 
 \smallskip
 
   We may thus apply to ${\bf Ref}_R$ the results proved by [QSSS 76], [QSS 79] and [K 91 Chap.II, \S6]. The main one, for our purposes, is the {\it cancellation theorem}:
   
   \medskip

\n {\bf Theorem 2.5.2} - {\it  Let $(M,h), \  (M',h')$ and
$(N,g)$ be $\epsilon$--hermitian spaces such that $:$  $$(M,h) \oplus (N,g) \ \simeq \ (M',h') \oplus (N,g).$$

\n Then
$(M,h) \simeq (M',h')$.}  

\n [Of course it is assumed that  $R$ and the $R$--modules $M,M',N$  
satisfy the conditions in $\S2.2$, so that $M,M'$ and
$N$ belong to the category ${\bf Ref}_R$.]

\smallskip
\n {\it Proof.} This is a special case of  theorem 6.6.1 of [K 91, chap.II]. The proof given there has two parts:

 a)  when $M,M'$ and $N$ are isotypic
of the same self-dual type, by reduction to the standard Witt
cancellation theorem for hermitian spaces over division algebras, cf. e.g. [A9, \S4.3, th.1] or [K 91, chap.I, th.6.5.2].

 b) in the general case, by using the isotypic decompositions of hermitian spaces, cf. [K 91, chap.II, th.6.3.1].
   
  \medskip

   {\bf 2.6.  Witt groups and Grothendieck-Witt groups}
   
   \medskip
   
   The {\it Grothendieck-Witt group}  ${\rm WGr}^\epsilon (R)$ is defined as
   the Grothendieck group of the category of the $\epsilon$--hermitian spaces over  $R$ (having property (2.2.3), of course), with respect to the orthogonal direct sum. If $(M,h)$ is an $\epsilon$--hermitian space, its class in ${\rm WGr}^\epsilon (R)$ will be denoted by $[M,h]_{gr}$. Because
   of theorem 2.5.2, we have:
   
   \smallskip
   
 \ (2.6.1) \quad    $ (M,h) \simeq (M',h')  \ \Longleftrightarrow \ [M,h]_{gr} = [M',h']_{gr} \ .$

   \smallskip
   
   Let ${\rm Hyp}^\epsilon_R$ be the subgroup of ${\rm WGr}^\epsilon (R)$ generated
   by the classes of the hyperbolic spaces (cf. \S2.4). The quotient
   
   \smallskip
   
  \quad ${\rm W}^\epsilon (R) \ = \ {\rm WGr}^\epsilon (R)/{\rm Hyp}^\epsilon_R$
   
   \smallskip

\n   is called the {\it Witt group} of $(R,\epsilon)$. If $(M,h)$  is an $\epsilon$--hermitian space, its class in ${\rm W}^\epsilon (R)$ will be denoted by  $[M,h]$.

\medskip

\n {\bf Proposition 2.6.1} - {\it Let $(M,h)$  and $(M',h')$  be two $\epsilon$--hermitian spaces. If $M \simeq M'$ and $[M,h] = [M',h']$, then $(M,h) \simeq (M',h')$.}

\smallskip

\n {\it Proof.} Since $[M,h] = [M',h']$, there exist two hyperbolic  $\epsilon$--hermitian spaces $(N,q)$ and $(N',q')$ such that 
   
   \smallskip
   
  \quad  $(M,h) \oplus (N,q) \ \simeq \ (M',h') \oplus (N',q')$.
   
   \smallskip
   
\n   This implies $M \oplus N \simeq M' \oplus N'$, and since $M \simeq M'$,
   we have $N \simeq N'$; by prop.2.4.1, we have $(N,q) \simeq (N',q')$,
   hence $(M,h) \simeq (M',h')$ by theorem 2.5.2.
   
   \medskip
   
   \n {\bf Proposition 2.6.2} - (i) {\it The subgroup ${\rm Hyp}^\epsilon_R$  of ${\rm WGr}^\epsilon (R)$ is torsion-free.}
   
   \smallskip
  (ii) {\it Let  $\pi:  {\rm WGr}^\epsilon (R) \ \to \ {\rm Hyp}^\epsilon_R$ be the homomorphism which maps the class of $(M,h)$ to the class of the
  hyperbolic space $H_\epsilon(M)$. Then $\pi(x) = 2x$ for every $x\in {\rm Hyp}^\epsilon_R.$}
   \smallskip
   
   \n {\it Proof of} (i). Let $x \in {\rm Hyp}^\epsilon_R $ be such that $nx = 0$,
   with $n > 0$. There exist reflexive modules $N$ and $N'$ such that
   $x = [H(N)]_{gr} - [H(N')]_{gr}$; the relation $nx=0$ means that $H(n\bullet N) \simeq H(n\bullet N')$, where $n\bullet N$ (resp. $n\bullet N'$) denotes the direct sum of $n$ copies of $N$ (resp. of $N'$). This implies
 $n\bullet (N\oplus N^*) \simeq n\bullet (N'\oplus N'^*)$; by the Krull-Remak-Schmidt theory for modules, we have $N \oplus N^* \simeq N' \oplus N'^*$, hence
 $H(N) \simeq H(N')$ by proposition 2.4.1; hence $x=0$.
 \smallskip
 
 \n {\it Proof of} (ii). It is enough to prove the formula when $x$ is the class
 of an hyperbolic space $H_\epsilon(N)$; in that case $\pi(x)$ is the class of $H_\epsilon(N \oplus N^*) = H_\epsilon(N) \oplus H_\epsilon(N)$, which is $2x$.
 \smallskip
    
   \

   \n  {\bf Proposition 2.6.3} - {\it Let $n$ be an integer $>0$. The following properties are equivalent} :
   
   \smallskip
   
\n \ \ (a) {\it ${\rm WGr}^\epsilon (R)$ is $n$--torsion free.}

    \smallskip
 \n \ \ (b) {\it If $(M,h)$ and $(M',h')$ are two $\epsilon$--hermitian spaces such that $n\bullet(M,h) \simeq n\bullet (M',h')$, then $(M,h) \simeq (M',h')$.
 
   \smallskip
  If $n$ is odd, these properties are also equivalent to} :
  
   \smallskip
   \n \ \  (c) {\it ${\rm W}^\epsilon (R)$ is $n$--torsion free}. 
   
  \smallskip
   \n   [Recall that an abelian group  $A$  is called {\it $n$--torsion free} if
 $a\in A$ and $na=0$ imply $a=0$.]
 
     \smallskip
     
     \n {\it Proof of} (a) $\Longleftrightarrow$ (b). This follows from (2.6.1).
     
     \smallskip
     
     \n {\it Proof of} (a) $\Longrightarrow$ (c). We have to show that, if $x \in  {\rm WGr}^\epsilon (R)$ and $nx \in  {\rm Hyp}^\epsilon_R$, then $x$ belongs to ${\rm Hyp}^\epsilon_R$. By part (ii) of proposition 2.6.2 we have $\pi(nx)=2nx$, hence $n(2x - \pi(x))=0$. Since ${\rm WGr}^\epsilon (R)$ is $n$--torsion free, this shows that $2x = \pi(x)$. Hence $2x$ belongs to ${\rm Hyp}^\epsilon_R$; since
     $nx$ has the same property, and $n$ is odd, we have $x\in {\rm Hyp}^\epsilon_R$.
     
     \smallskip 
     
      \n {\it Proof of} (c) $\Longrightarrow$ (a). Let $x 
     \in {\rm WGr}^\epsilon (R)$ be such that $nx = 0$. Since ${\rm W}^\epsilon (R)$ is
     $n$--torsion free, we have $x \in {\rm Hyp}^\epsilon_R$, and, by part (i) of proposition 2.6.2, we have $x=0$.
     
     \smallskip

   In the next \S, we shall see that {\it properties} (a), (b) {\it and} (c) {\it hold for every odd $n$ and for every $(R,\epsilon)$.}
        \medskip
     
\centerline
{\bf \S 3. Division properties for $\epsilon$--hermitian spaces}

\medskip

We keep the notation of \S 2. In particular, $k$ is a field of characteristic $\not = 2$,   $R$ is a finite dimensional $k$--algebra
with a $k$--linear involution $x \mapsto \overline x$, and $\epsilon = \pm 1$.

\medskip

{\bf 3.1. Statement of the division theorem }

\smallskip

Let $(M,h)$ be an $\epsilon$--hermitian space over $R$, and let $n$ be an  integer $>0$. We say that $(M,h)$ {\it has property} ${\rm Div}_n$ if
every $(M',h')$ with $n\bullet (M',h') \simeq n\bullet (M,h)$ is isomorphic to $(M,h)$, i.e. if one is allowed to `` divide by  $n$  ''. [Recall that $n\bullet (M,h)$ is the direct sum of $n$ copies of $(M,h)$.] Note that
the hypothesis $n\bullet (M',h') \simeq n\bullet (M,h)$ implies $M \simeq M'$
by the Krull-Remak-Schmidt theorem for modules; hence we may assume
that $M'=M$, so that the problem consists in comparing two different $\epsilon$-hermitian forms on the same $R$-module.

\smallskip

We say that $R$ {\it has property} ${\rm Div}_{n,\epsilon}$ if every $\epsilon$--hermitian space over $R$ has property ${\rm Div}_n$.  The goal of $\S3$ is to prove:

\smallskip

\n {\bf Theorem 3.1.1} - {\it Let $n$ be a positive odd integer. Then} :

 (a) {\it The group ${\rm WGr}^\epsilon (R)$ is $n$--torsion free.}
 
 (b) {\it The group ${\rm W}^\epsilon (R)$ is $n$--torsion free}.
 
 (c) {\it Property ${\rm Div}_{n,\epsilon}$ holds}.
 
 \smallskip
 
 By proposition 2.6.3 the three assertions (a), (b) and (c) are equivalent. Their proof will consist in a series of reductions.
  
  \smallskip

     \n {\it Remark.} Division by 2 is rarely possible. For instance,
     if $k = \Q$, the quadratic forms  $ x^2$ and $2x^2$ are not isomorphic,
     but they become so after duplication:
     
       $x^2 + y^2 \ \simeq \ (x+y)^2 + (x-y)^2 = 2x^2 + 2y^2 .$
       
       \n However, we shall see later (theorem 3.5.3) that, if $k$ is real closed, then ${\rm WGr}^\epsilon (R)$ is torsion-free, so that property ${\rm Div}_{n,\epsilon}$ holds for every  $n > 0$, including $n=2$.

  \medskip
  
  {\bf 3.2. The ${\rm W}(k)$--module structure of ${\rm W}^\epsilon(R)$}
  
  \smallskip
  
  Let $(M,h)$ be an $\epsilon$--hermitian space over $R$, and let $V$
  be a finite-dimensional $k$--vector space equipped with a nonsingular symmetric bilinear form $q$ (note that $(V,q)$ may be viewed as a 1--hermitian space over $k$). The tensor product  $ V \otimes_k M$ has a
  natural structure of $\epsilon$--hermitian space over $R$, which we denote by $(V,q) \otimes_k (M,h)$: this is easily checked, using the fact that the $R$-dual of $ V \otimes_k M$ is canonically isomorphic to $ V' \otimes_k M^*$, where $V'$ is the $k$--dual of $V$.
  
  \smallskip

  Let ${\rm WGr}(k)$ denote the Witt-Grothendieck ring of the field $k$. Using the construction above, we see that ${\rm WGr}^\epsilon(R)$ {\it has a natural structure of ${\rm WGr}(k)$-module}. After dividing by hyperbolic spaces, this shows that ${\rm W}^\epsilon(R)$ {\it has a natural structure of ${\rm W}(k)$-module}.
We shall use this structure to prove:

\medskip

\n {\bf Proposition 3.2.1} - {\it Suppose that  $k$ is not formally real. Then} :

 (a) {\it There exists $r \in \N$  such that $2^r.{\rm W}^\epsilon(R)= 0$.} 
 
 (b) {\it Property ${\rm Div}_{n,\epsilon}$ holds for every odd $n$ and every $\epsilon = \pm1$.}

\smallskip
\n[Recall that  $k$ is  {\it formally real} if $-1$ is not a sum of squares
in $k$; this is equivalent to the existence of a total order on $k$ which is compatible with the field structure, cf.  [A6, \S2.2]; such a field has characteristic 0.]

\smallskip

\n {\it Proof.} We use the same method as in [La 04, p.197 and p.255]. The hypothesis on $k$ implies that there exists $r \in \N$ such that $2^r = 0$ in ${\rm W}(k)$, cf. e.g. [S 85, chap.II, th.6.4 (ii)]. Since ${\rm W}^\epsilon(R)$ is a ${\rm W}(k)$--module, we have  $2^rx=0$ for every $x \in {\rm W}^\epsilon(R)$. This proves (a). As for (b), it follows from (a) and from proposition 2.6.3, since ${\rm W}^\epsilon(R)$ is $n$-torsion free.

\smallskip

Proposition 3.2.1 allows a first reduction: {\it it is enough to prove theorem 3.1.1 when $k$ is formally real}. Before we do the second reduction, we need
to recall some elementary (but not completely obvious) properties of hermitian spaces with respect to field extensions.

\medskip

{\bf 3.3. Field extensions : going up }

\smallskip

Let $K$ be a field extension of $k$ and let $R_K = K \otimes_k R$. If $M$
is an $R$-module, let $M_K$ denote the $R_K$-module $K \otimes_k M$.
 Let $M^*={\rm Hom}_R(M,R)$ be the dual
of $M$ and let $(M_K)^*= {\rm Hom}_{R_K}(M_K,R_K)$ be the dual of $M_K$. Every element of $M^*$ extends lineary to an element of $(M_K)^*$ and we obtain in this way a natural map $$\iota_M : K \otimes M^* \to (M_K)^*.$$

\n {\bf Proposition 3.3.1} - {\it If $M$ is a finitely generated $R$-module,  $\iota_M$ is an isomorphism.} 

\smallskip

\n {\it Proof}. Let $a$ and $b$ denote the functors $M \mapsto K \otimes M^*$ and $M \mapsto (M_K)^*$. When $M = R$, we have $a(M) = b(M)=R_K$
and $\iota_M$ is the identity. By additivity, this shows that the proposition is true when $M$ is $R$-free. In the general case, we choose an exact 
sequence  $L' \to L \to M \to 0$, where  $L'$ and $L$ are finitely generated $R$-free modules; this is possible since $M$ is finitely generated and the ring $R$ is noetherian. Since both $a$ and $b$ are left exact functors, we have a commutative diagram where the vertical maps are  $0, \iota_M, \iota_L$ and $\iota_{L'} $ \ :

\begin{center}
 \ $\ 0 \ \rightarrow \ a(M)  \ \rightarrow  \ a(L)   \  \rightarrow  \ a(L')$
   
$\!    \downarrow  \ \ \ \ \qquad \downarrow \ \  \ \quad \qquad \downarrow  \ \ \  \quad \quad \ \ \downarrow  $

 $\ \ \ 0 \ \rightarrow \ b(M) \  \rightarrow \  b(L) \ \ \rightarrow  \ b(L') . $

\end{center}
\smallskip

\n Since $\iota_L$ and $\iota_{L'}$ are isomorphisms, so is $\iota_M$.

\medskip

  The proposition above allows us to identify $K \otimes M^*$ with $(M_K)^*$
  and to denote both modules by  $M_K^*$.
  Let now $(M,h) $ be an $\epsilon$-hermitian space over $R$. The isomorphism $h : M \to M^*$ gives, by extension of scalars, an isomorphism $h_K : M_K \to M_K^*$. We thus obtain an $\epsilon$-hermitian space $(M_K,h_K)$ over $R_K$. The functor $(M,h) \mapsto (M_K,h_K)$ transforms hyperbolic spaces into hyperbolic spaces. Hence we get natural maps $$ {\rm WGr}^\epsilon(R) \to {\rm WGr}^\epsilon(R_K)   \ \ {\rm and} \ \ {\rm W}^\epsilon(R) \to {\rm W}^\epsilon(R_K).$$
  
\n  In [S 85, chap.2, \S5], these maps are denoted by $r_{K/k}^*$. We shall
 simplify the notation by writing them $x \mapsto x_K$.
 
    \medskip  

{\bf 3.4. Field extensions : going down}

\medskip
In this section, $K$ is a finite field extension of  $k$. We choose a non-zero
$k$-linear map $s: K \to k$.

\smallskip
\n {\it The Scharlau transfer for vector spaces}

\smallskip
  Let $M$ be a vector space over $K$, and let $N$ be a finite dimensional vector space over $k$. Let $N_K = K \otimes_k N$; the map $s$ defines $s_N = s\otimes1 : N_K \to N$. If $f: M \to N_K$ is $K$-linear, the composition $s_N \circ f$ is a
  $k$-linear map $M \to N$. We thus get a map $$s_*:  {\rm Hom}_K(M,N_K) \to {\rm Hom}_k(M,N).$$

  \smallskip
  \n {\bf Proposition 3.4.1} - {\it The map $s_*: {\rm Hom}_K(M,N_K) \to {\rm Hom}_k(M,N)$ is an isomorphism.}
  
\smallskip

\n {\it Proof.} Put $a(M,N)={\rm Hom}_K(M,N_K)$ and $b(M,N)={\rm Hom}_k(M,N)$. Both functors transform direct sums (for $M$) into direct products, and direct sums (for $N$) into direct sums. Hence it is enough
to check that $s_*$ is an isomorphism when $M$ and $N$ are 1-dimensional, i.e. when  $M = K$ and $N=k$; in that case, both $a(M,N)$
and $b(M,N)$ are $K$-vector spaces of dimension 1, and the map is $K$-linear and non zero; hence it is an isomorphism. 

\n {\it Remark.} When $K = \C, k=\R$ and $N=k$, the proposition amounts to the standard fact that every $\R$-linear form on a complex vector space is the real part of a unique $\C$-linear form.

\smallskip

\n {\it The Scharlau transfer for modules}
\smallskip

  We keep the same notation $(K,s,N,M)$ as above, and we assume
  that both $N$ and $M$ are modules over a $k$-algebra $A$; in that case,
  $M$ is an $A_K$-module, where $A_K = K \otimes_k A$.  We have natural inclusions:
  \smallskip
  $ {\rm Hom}_{A_K}(M,N_K) \subset {\rm Hom}_K(M,N_K) \quad {\rm and} \quad
  {\rm Hom}_{A_k}(M,N) \subset {\rm Hom}_k(M,N).$
  
  \n The map $s_*$ of proposition 3.4.1 maps the first subspace into the
  second one. Hence we get
    $ s_{*,A}: {\rm Hom}_{A_K}(M,N_K) \to {\rm Hom}_A(M,N).$

  \smallskip
  
  \smallskip
  \n {\bf Proposition 3.4.2} - {\it The map $s_{*,A}$ is an isomorphism.}
\smallskip

\n {\it Proof.} The injectivity of $s_{*,A}$ follows from that of $s_{*}$. For the surjectivity, it is enough to show that, if $f: M \to N_K$ is a $K$-linear map such that $s_N \circ f$ is $A$-linear, then so is $f$. Indeed, if $a$ is any element of  $A$, the two maps  $x \mapsto f(ax)$ and $x \mapsto af(x)$ have the same image by $s_*$, hence they coincide.

\smallskip

\smallskip

\n {\it The Scharlau transfer for hermitian spaces}
\smallskip

We apply the above to $A=R$. If $(M,h)$ is an $\epsilon$-hermitian space
over $R_K$, the $K$-bilinear map  $h_{{\rm b}} : M \times M \to R_K$,
composed with  $s_R : R_K \to R$, gives an $\epsilon$-hermitian form
$s_*\circ h_{{\rm b}} $ on  $M$, where $M$ is now viewed as an $R$-module.

\smallskip
\n {\bf Proposition 3.4.3} - {\it The hermitian form $s_*\circ h$ is nonsingular.}

\smallskip

\n {\it Proof.} We have to show that the corresponding map $M \to {\rm Hom}_R(M,R)$ is an isomorphism. This map factors in:

\smallskip

  \quad  $ M \ \to \ {\rm Hom}_{R_K}(M,R_K) \ \to \ {\rm Hom}_R(M,R)$.

\smallskip

The left side map is an isomorphism because $h$ is nonsingular; the right
side map is an isomorphism because of proposition 3.4.2. Hence the
composite is an isomorphism.

\smallskip

  We thus obtain an $\epsilon$-hermitian space over $R$, which is called
  the {\it Scharlau transfer} of  $(M,h)$, relative to $s$. We denote it by
  $(M,s_*h)$. We use a similar notation for the corresponding homomorphisms
  
  $$ s_* :  {\rm WGr}^\epsilon(R_K) \to {\rm WGr}^\epsilon(R)\ \ \ {\rm and} \ \  s_* : {\rm W}^\epsilon(R_K) \to {\rm W}^\epsilon(R).$$
  
  \smallskip
  
  \n {\it A basic formula}
  
    \smallskip
  
  The definitions above apply to $R=k$ and they give
  a map $s_*: {\rm W}(K) \to {\rm W}(k)$ between the corresponding Witt groups. The element $s_*(1)$ of ${\rm W}(k)$ is the Witt class of the $k$-quadratic form  $z \mapsto s(z^2)$ on the $k$-vector space $K$.
  
   \smallskip
   
   \n {\bf Proposition 3.4.4} - {\it Let $x \in {\rm W}^{\epsilon}(R)$ and let
   $x_K$ be its image in ${\rm W}^{\epsilon}(R_K)$, cf.} \S 3.3. {\it We have} :
   
   \smallskip
 \qquad  $  s_*(x_K) = s_*(1).x  \ \ {\it in } \ \ {\rm W}^{\epsilon}(R).$
   
   \smallskip
   
   \n [The product of $s_*(1)$ and $x$ makes sense, because of the ${\rm W}(k)$-module structure of ${\rm W}^{\epsilon}(R)$, cf. \S3.2.]
   
   \smallskip
   \n {\it Proof.} It is enough to check the formula when $x$ is the class of
   an $\epsilon$-hermitian space $(M,h)$. In that case, the $k$-bilinear map
   $M_K \times M_K \to R$ corresponding to $s_*(M_K,h_K)$ is given by:
   
   \smallskip
   
     $(z \otimes m) \times (z' \otimes m') \ \mapsto \ s(zz')h(m,m')$  \quad  for 
     $z,z' \in K$  and  $m,m' \in M$.
     
     \n This is the same formula as for $(K,q_s) \otimes (M,h)$, where $q_s$
     is the quadratic form  $z \mapsto s(z^2)$ on $K$.

   \medskip
 
\n {\it An application of the basic formula}

\medskip
\n {\bf Proposition 3.4.5} - {\it Let $x \in {\rm W}^{\epsilon}(R)$.}

(i) {\it Let \ $a \in k^\times$ 
 be such that neither \ $a$ \ nor \ $-a$\ is a square and let \ $K = k(\sqrt a), K' = k(\sqrt{-a})$. If $x_K= 0$ and $x_{K'}=0$, then $2x=0$.}

(ii) {\it Let $E$ be a finite extension of $k$ of odd degree. If $x_E=0$, then $x=0$.}

\smallskip

\n {\it Proof of}  (i). Let $s: K \to k$ be the trace map divided by 2. The corresponding quadratic form $s_*(1)$ is the form $X^2 + aY^2$, i.e. $<\!1\!> + <\!a\!>$ using standard notation. Similarly, there is a linear map $s' : K' \to k$ such that
$s'_*(1) = \ <\!1\!> + <\!-a\!>$; we have $$s_*(1)+s'_*(1)\ = \ <\!1\!> + <\!1\!> + <\!a\!> + <\!-a\!> \ = \ <\!1\!> + <\!1\!> \ =  \ 2,$$ because $<\!a\!> + <\!-a\!> \ =0$  in ${\rm W}(k)$. By proposition 3.4.4, the relations $x_K=0$ and $x_{K'}=0$
imply that $s_*(1)x=0$ and $s'_*(1)x=0$; by adding up, we find $2x=0$.

\smallskip

\n {\it Proof of} (ii). Since $[E:k]$ is odd, a well-known construction of Scharlau (see e.g. [S 85, chap.2, lemma 5.8]) shows that there exists a linear form $s:E \to k$ such that $s_*(1) = 1$ in ${\rm W}(k)$.
The relation $x_E = 0$ implies $s_*(1) x = 0$, i.e. $x=0$.

\smallskip

\n {\it Remark.} Part (ii) of proposition 3.4.5 is equivalent to the main result of [BFL 90].

\medskip
{\bf 3.5. Second reduction: from formally real to real closed}

\medskip
\n
{\bf Proposition  3.5.1} - {\it Let $\Omega$ be an
algebraic closure of $k$. Let $x \in {\rm W}^{\epsilon}(R)$. Suppose that $x_K= 0$ for every
real closed extension $K$ of $k$ contained in $\Omega$. Then $x$ is a $2$-primary torsion element, i.e.
there exists $r \in {\bf N}$ such that $2^r  x = 0$}. 

\n [Recall that a field is {\it real closed} if it is formally real, and if no proper algebraic extension of it has this property; such a field is called {\it ordonn\'e maximal} in [A6, \S2].]

\medskip
\n
{\it Proof.} Suppose that no such $r$ exists. By Zorn's lemma, we may choose an extension $L \subset \Omega$ such that $x_L$ is not a 2--primary torsion element, and such
that $L$ is maximal for this property. Then  prop.3.2.1 implies that $L$ is formally real. Let us show that $L$ is real closed, which will give a contradiction. To do so it is enough to prove (cf. [A6, \S2, th.3]) that $L$
has the following properties:

 (i) Every element  $a$ of  $L$ is either a square or minus a square.
 
 (ii) The field $L$ does not have any non trivial odd degree extension.
 
 Let us prove (i). If not, there would exist $a\in L$ such that neither $a$
 nor $-a$ is a square. The field extensions $L( \sqrt a)/L$ and $L(\sqrt {-  a})/L$ are non trivial. Because of the maximality of $L$, this implies that 
$x_{L( \sqrt a)}$ and $x_{L( \sqrt {- a})}$ are 2--primary torsion elements. 
By part (i) of proposition 3.4.5 the same is true for $x_L$, which is a contradiction.

The proof of (ii) is similar: one uses part (ii) of proposition 3.4.5.


\smallskip

\n {\bf Corollary 3.5.2} - {\it Let $n$ be an odd integer. If $R_K$ has property ${\rm Div}_{n,\epsilon}$ for every real closed field extension $K$ of $k$, then $R$ has the same property.}

\smallskip

\n {\it Proof.} Let $x \in {\rm W}^\epsilon(R)$ be such that $nx = 0$. If $K$ is any real closed extension of $k$, we have $x_K$ = 0, because of ${\rm Div}_{n,\epsilon}$. By proposition 3.5.1 this implies that  $2^rx = 0$ for a suitable $r\in \N$. Since $nx=0$ and $n$ is odd, we have $x=0$.

\smallskip

Corollary 3.5.2 shows that {\it it is enough to prove theorem 3.1.1 when $k$ is real closed}.  We shall actually prove more:

\smallskip

\n {\bf Theorem 3.5.3} - {\it If $k$ is real closed, then $R$ has property ${\rm Div}_{n,\epsilon} $ for every $n > 0$ and every $\epsilon = \pm 1.$}

\smallskip

  In the next section, we shall check this when $R$ is a division algebra. The general case will be done in \S3.11.

\medskip

{\bf 3.6. The case where $k$ is real closed and $R$ is a division algebra}

\smallskip
\n {\bf Theorem 3.6.1} - {\it Suppose that $k$ is real closed and that $R$ is a division algebra.  Then ${\rm WGr}^\epsilon(R)$ is isomorphic to either $\Z$
or $\Z \oplus \Z$.}

\smallskip

  In particular, ${\rm WGr}^\epsilon(R)$ is torsion-free. By proposition 2.6.3 this implies:
  
  \smallskip
  
  \n {\bf Corollary 3.6.2} - {\it Property ${\rm Div}_{n,\epsilon} $ is true for every $n > 0$ and every $\epsilon = \pm1$.}
  
  \smallskip
  
  \n {\it Proof of theorem 3.6.1.} Up to isomorphism, $R$ is either 
    $k,  k(i)$ or the standard quaternion algebra over $k$. There are ten different possibilities, according to the choice of $R$, of its involution, and of  $\epsilon$ :
 \smallskip
 
  3.6.3 - $R = k$ with trivial involution and $\epsilon = 1$. Then  ${\rm WGr}^\epsilon(R)$ is the 
  Witt-Grothendieck group of quadratic forms over $k$. Such a form is well-defined by its signature, i.e. by the number of positive and negative 
  coefficients when written as  ${\Sigma} a_ix_i^2$. Hence ${\rm WGr}^\epsilon(R) \simeq \Z \oplus \Z.$
  \smallskip

  3.6.4 - $R = k$ with trivial involution and $ \epsilon = -1$. The classification of alternating forms shows that such a form is well-defined by its rank, which is an even number. Hence
 ${\rm WGr}^\epsilon(R) \simeq \Z.$

\smallskip

3.6.5 - $ R = k(i)$ with trivial involution and $\epsilon = 1$. Here ${\rm WGr}^\epsilon(R) $ is the same as the Witt-Grothendick group of quadratic forms over $k(i)$;
such a form is well-defined by its rank, hence ${\rm WGr}^\epsilon(R)  \simeq \Z.$

\smallskip

3.6.6 - $ R = k(i)$ with trivial involution and $\epsilon = -1$. The result is the same as in case 3.6.4, i.e.  ${\rm WGr}^\epsilon(R)  \simeq \Z.$

\smallskip

3.6.7 - $ R = k(i)$ with the involution $i \mapsto -i$ and $ \epsilon = 1$.
 The result is the same as in case 3.6.3, i.e. ${\rm WGr}^\epsilon(R) \simeq \Z \oplus \Z.$

\smallskip

3.6.8 - $ R = k(i)$ with the involution $i \mapsto -i$ and $ \epsilon = -1$. If $h$ is an $\epsilon$-hermitian form, then $i.h$ is $(-
\epsilon)$-hermitian. We are thus reduced to case 3.6.7, hence ${\rm WGr}^\epsilon(R) \simeq \Z \oplus \Z.$

\smallskip

3.6.9 - $R = $ quaternion algebra with basis $<\!\!1,i,j,ij \! \!>$, the involution being the hyperbolic one $\{i,j,ij \mapsto  -i,-j,-ij\}$, and $\epsilon = 1$. The result is the same as in cases 3.6.3 and 3.6.7, namely ${\rm WGr}^\epsilon(R) \simeq \Z \oplus \Z.$

\smallskip
3.6.10 - $R $ = quaternion algebra with the hyperbolic involution and $\epsilon = -1$. Up to isomorphism, there is only one form of rank 1,  namely
$(x,y) \mapsto xi\overline y$. Hence ${\rm WGr}^\epsilon(R)  \simeq \Z.$

\smallskip

3.6.11 - $R$ = quaternion algebra with an orthogonal involution, for instance $\{i,j,ij \mapsto  -i,j,ij\}$, and $\epsilon = 1$.  If $h$ is an $\epsilon$-hermitian form, then $hi$ is $(-\epsilon)$-hermitian for the standard involution of $R$. By 3.6.10, we thus have ${\rm WGr}^\epsilon(R)  \simeq ~\Z.$

\smallskip

3.6.12 - $R$ = quaternion algebra with an orthogonal involution as above and $\epsilon = -1$.  The same argument as in 3.6.11 shows that ${\rm WGr}^\epsilon(R) \simeq \Z \oplus \Z.$

\smallskip

This completes the proof of theorem 3.6.1.

\smallskip

\n {\it Remark.} The corresponding Witt groups ${\rm W}^\epsilon(R)$ are easily computed. They are respectively: 
\quad $\Z, \ 0, \ \Z/2\Z, \ 0, \ \Z, \ \Z,\ \Z, \ \Z/2\Z,\ \Z/2\Z,\ \Z.$

\medskip

{\bf 3.7. Hermitian elements}

\smallskip

Before completing the proof of theorems 3.1.1 and 3.5.3, we recall a few basic facts on hermitian elements.

\smallskip
 Let $E$ be a ring with an involution $\sigma : E \to E$ \footnote{In the next section,
  $E$ will not be $R$, but rather the endomorphism ring of an hermitian $R$-space; in Bourbaki's parlance ([A8, \S1.3]), we switch from a module to its ``contre-module"; from a ``category with duality" point of view, we switch from an object to its endomorphism ring.}. Let $\epsilon = \pm 1$, and put 
$$E^{\epsilon} = \{ z \in E^{\times} \ | \ \sigma(z) = \epsilon z \}.$$  
If $z \in E^{\epsilon}$, the map $h_z : E \times E \to E$ defined by
$h_z(x,y) = x.z.\sigma(y)$ is an $\epsilon$--hermitian space over $E$;
conversely, every $\epsilon$--hermitian space over $E$ with underlying module $E$ is isomorphic to $h_z$
for some $z \in E^{\epsilon}$.

Define an equivalence
relation on $E^{\epsilon} $ by setting 
$z \equiv z'$ if there exists $e \in E^\times$ with $z' = \sigma(e) z e$; this is equivalent to $ (E,h_z) \simeq (E,h_{z'})$. Let
$ H^{\epsilon}(E,\sigma)$
\n be the quotient of $E^{\epsilon}$ by this equivalence relation.
If $z \in E^{\epsilon} $, we denote by $[z]$ its class in $H^{\epsilon}(E,\sigma)$. 

\medskip
\n
{\bf Lemma 3.7.1} - {\it Let $A$ be a ring, let $A^{{\rm op}}$ be the opposite ring, and let $E = A \times A^{{\rm op}}$. Let $\sigma$ be the involution of $E$ defined by $(a,b) \mapsto (b,a)$.
Then $H^{\epsilon}(E,\sigma) = \{1 \}$.}

\smallskip
\n
{\it Proof.} Every element of $E^\epsilon$ is of the form $(a,\epsilon a)$ with 
$a \in A^\times$. If $(b,\epsilon b)$ is another such element, we have $(a^{-1},b) (a,\epsilon a) (b, a^{-1}) = (b,
\epsilon b)$,
hence $(a, \epsilon a) \equiv (b, \epsilon b)$. 

\medskip
\n
{\bf Lemma 3.7.2} - {\it  Let $(A,\tau)$ and $(B,\tau')$ be two rings with
involution, and let $(E,\sigma)$ be their direct product.   Then $H^{\epsilon}(E,\sigma) = H^{\epsilon}(A,\tau) \times H^{\epsilon}(B,\tau')$. }

\smallskip
\n
{\it Proof.} Let $(a,b), (a',b') \in E^{\epsilon} $. Then $[(a,b)] = [(a',b')]$ in $H^{\epsilon}(E,\sigma)$
if and only if there exists $(e,f) \in E^{\times}$ such that $\sigma(e,f)(a,b)(e,f) = (a',b')$.
This is equivalent to $\tau (e) a e = a'$ and $\tau'(f) b f = b'$, hence to
$[a] = [a'] \in H^{\epsilon}(A,\tau)$ and $[b] = [b'] \in H^{\epsilon}(B,\tau')$. The lemma follows.  

\medskip
\n
{\bf Lemma 3.7.3} - {\it Let $J$ be a two--sided ideal  of $E$ such that $\sigma(J) = J.$ Suppose that $J$ is
nilpotent, and that $2$ is invertible in $E$. Set $\overline E = E/J$ and denote by $\overline \sigma$ the involution of  $\overline E$ defined by $\sigma$. Then the projection map $ \pi : E \to \overline E$
induces a bijection \ 
$H^{\epsilon}(E,\sigma) \to H^{\epsilon}(\overline E,\overline \sigma).$}

\smallskip
\n
{\it Proof.}  Let us first show the surjectivity of the projection $E^{\epsilon} \to \overline E^{\epsilon}$ and hence that of $H^{\epsilon}(E,\sigma) \to H^{\epsilon}(\overline E,\overline \sigma).$
Let $u$ be an element of $\overline E^{\epsilon}$. Choose $v \in E$ with $\pi(v)=u$ and let $w={1\over2}(v + \epsilon \sigma(v))$. We have
$\sigma(w) = \epsilon w$ and $\pi(w)=u$. Since the kernel $J$ of $\pi$ is nilpotent, and $u$ is invertible, the relation $\pi(w) = u$ implies that $w$ is invertible. Hence we have $w \in E^{\epsilon}$; this proves the surjectivity of $E^{\epsilon} \to \overline E^{\epsilon}$.

Let us now prove the injectivity of $H^{\epsilon}(E,\sigma) \to H^{\epsilon}(\overline E,\overline \sigma).$ By induction, we may assume that $J^2 =  0$. Let $z, z' \in E^{\epsilon}$
such that $[z] = [z']$ in $H^{\epsilon}(\overline E,\sigma).$ Then there exists
$a \in E^{\times}$ such that $\sigma(a) z a - z' \in J$. Replacing $z$ by $\sigma(a) z a $,
we may assume that $z' = z + r$ with $r \in J$. Then we have $\sigma(r) = \epsilon r$.
Set $b = {1 \over 2} z^{-1} r$. Then  $z b = {1 \over 2}r$, hence $\sigma(b) z = {1 \over 2}r$; by adding up, this gives
 $r = zb + \sigma(b) z$. Since $J^2=0$, we have $\sigma(b) z b = 0$, hence 
$z' = \sigma(1 + b) z (1 + b)$; therefore $[z] = [z']$ in $H^{\epsilon}(E,\sigma) $,
as claimed.

\medskip

{\bf 3.8. Classifying hermitian spaces via hermitian elements}

\medskip
Let $\epsilon_0 = \pm 1$, and let $(M,h_0)$ be 
an $\epsilon_0$--hermitian space over $R$. Set $E_M = {\rm End}(M)$. Let $\tau : E_M \to E_M$ be the involution of $E_M$ {\it induced} by $h_0$, i.e. 

$$\tau (e) = h_0^{-1} e^* h_0,  \quad {\rm for} \ \  e \in E_M, $$

\n where $e^*$ is the adjoint of  $e$, cf. \S2.1. If $(M,h)$ is an $\epsilon$--hermitian space (with the same underlying module $M$), we have

\smallskip

\quad \quad   $\tau(h_0^{-1} h) = h_0^{-1}(h_0^{-1}h)^*h_0 = h_0^{-1}h^*(h_0^{-1})^*h_0 = \epsilon_0 \epsilon h_0^{-1} h$.

\smallskip

\n Hence $h_0^{-1} h$ is an $\epsilon_0 \epsilon$--hermitian element of $(E_M,\tau)$; let $[h_0^{-1} h]$ be its class in
$H^{\epsilon_0 \epsilon}(E_M,\tau).$

\medskip
\n
{\bf Lemma 3.8.1} - {\it Sending an $\epsilon$--hermitian space $(M,h)$ to the element $[h_0^{-1} h]$ of
$H^{\epsilon_0 \epsilon}(E_M,\tau)$ induces a bijection between the set of isomorphism classes
of $\epsilon$--hermitian spaces $(M,h)$ and the set $H^{\epsilon_0 \epsilon}(E_M,\tau)$. }

\smallskip
\n
{\it Proof.} Let $(M,h)$ and $(M,h')$ be two $\epsilon$--hermitian spaces, and set $u =  h_0^{-1}h, \ u' =h_0^{-1}h'$. We have $[u] = [u']$ if and only if there exists $e \in E_M^{\times}$ such that $u'= \tau(e)ue$, i.e. 
 
 \smallskip
 
\quad \quad   $h_0^{-1}h' = h_0^{-1}e^*h_0.h_0^{-1}h.e$ \
   
   \smallskip
   
   \n which is equivalent to $h' = e^*he$ and means that $(M,h)$ 	and $(M,h')$  are isomorphic.

\medskip

{\bf 3.9. A reformulation of the $n$-division property}

\smallskip
 Let $(M,h_0)$ be
an $\epsilon_0$--hermitian space, let $E = {\rm End}(M)$ and let $\sigma : E \to E$ be the involution induced
by $h_0$, cf. \S3.8. If $n \in {\bf N}$, let $E_n = {\rm End}(M^n) = {\bf M}_n(E)$, and
let $\sigma_n : E_n \to E_n$ be defined by $\sigma_n(a_{i,j})  = (\sigma(a_{j,i}))$.

\medskip
\n
{\bf Proposition 3.9.1} - {\it Let $n$ be an integer $> 0$.
The following properties are equivalent} :

\smallskip
\n
(i) {\it The $\epsilon$-hermitian spaces $(M,h)$ with underlying module $M$ have property ${\rm Div}_n$ of $\S3.1$. }
\smallskip

\n
(ii) {\it The map  $f_n :  H^{\epsilon_0 \epsilon}(E,\sigma) \to H^{\epsilon_0 \epsilon }(E_n,\sigma_n)$ given by 

$$[u] \mapsto  [u_n] =  \left [ \left( \matrix {u & 0 & \dots & 0 \cr \dots & \dots & \dots & \dots \cr
0  & \dots & 0 & u \cr } \right)  \right ] $$
\smallskip
\n
is injective, where $u_n$ is  the $n \times n$--matrix with diagonal entries equal
to $u$. }

\n [For the definition of $H^{\pm1}(E,\sigma)$, see \S3.7.]

\medskip
\n
{\it Proof.} Let us prove that (i) implies (ii). Let $u, u' \in E$ such that 
$f_n([u]) = f_n([u']).$ Let $(M,h)$, respectively $(M,h')$, be $\epsilon$--hermitian
forms such that the isomorphism class of $(M,h)$ corresponds to $[u] \in H^{\epsilon_0 \epsilon}(E,\sigma)$,
and that the isomorphism class of $(M,h')$ corresponds to $[u'] \in H^{\epsilon_0 \epsilon}(E,\sigma)$.
Then the isomorphism class of $n \bullet (M,h)$ corresponds to $f_n([u])$, and the
isomorphism class of $n \bullet (M,h')$  corresponds to $f_n([u'])$. 
By hypothesis, we have $f_n([u]) = f_n([u'])$, hence $n \bullet (M,h) \simeq n \bullet (M,h')$. By (i), we
have $(M,h) \simeq (M,h')$, thus $[u] = [u']$ in $H^{\epsilon_0 \epsilon}(E,\sigma)$. This implies that
$f_n$ is injective, hence (ii) holds. 

\smallskip
 Conversely, let us prove that (ii) implies (i). Let $(M,h)$ and $(M,h')$ be two 
$\epsilon$--hermitian spaces, and let $u, u' \in E$ such that $[u]$, respectively $[u']$,
are the elements of $H^{\epsilon_0 \epsilon}(E,\sigma)$ corresponding to the isomorphism
classes of $(M,h)$, respectively $(M,h')$, by the bijection of lemma 3.8.1. 
Then the isomorphism class of $n \bullet (M,h)$ corresponds to $f_n([u])$, and the
isomorphism class of $n \bullet (M,h')$ corresponds to $f_n([u']))$.
By hypothesis, we
have $n \bullet (M,h) \simeq n \bullet (M,h')$, hence $f_n(u) = f_n(u')$. By (ii), the map $f_n$ is
injective, therefore we have $(M,h) \simeq (M,h')$, and this implies (i).

\medskip

{\bf 3.10. An injectivity property}

\medskip
Let $E$ be a finite dimensional
$k$--algebra, and let $\sigma : E \to E$ be a $k$--linear involution. 
If $n \in {\bf N}$, let $E_n =  {\bf M}_n(E)$, and
let $\sigma_n : E_n \to E_n$ be defined by $\sigma_n(a_{i,j})  = (\sigma(a_{j,i}))$.

\medskip
\n
{\bf Theorem 3.10.1} - {\it Suppose that  $k$ is real closed and that $n$ is $>0$. Then the map \ $f_n :  H^{\epsilon}(E,\sigma) \to H^{\epsilon}(E_n,\sigma_n)$ \ given by 

$$[u] \mapsto  [u_n] =  \left [ \left( \matrix {u & 0 & \dots & 0 \cr \dots & \dots & \dots & \dots \cr
0  & \dots & 0 & u \cr } \right)  \right ] $$
\smallskip
\n
is injective, where $u_n$ is  the $n \times n$--matrix with diagonal entries equal
to $u$. }

\medskip
\n
{\it Proof.} 

\n {\it The case of a simple algebra.}  

 Suppose first that $E$ is a simple $k$--algebra. Then $E \simeq
{\bf M}_m(D)$ for some $m \in {\bf N}$ and some division algebra $D$.
Let $M = D^m$. Then
$E = {\rm End}_D(M)$.
It is well-known (see e.g. [S 85, chap.8, cor.8.3] or [KMRT 98, chap.I, th.3.1]) that there exists a $k$-linear involution of $D$ of the same kind as $\sigma$ (i.e. it is equal to $\sigma$ on the center of $D$). Let $\rho: D \to D$ be such an involution. By [KMRT 98, chap.I, th.4.2] there exist $\epsilon_0 = \pm 1$ and a nonsingular $\epsilon_0$-form $h_0$ on $M$ such that the
involution $\sigma$ of $E$ is induced by $h_0$ in the sense defined in \S3.8.

Let $u, u' \in E^{\epsilon}$ such that $f_n([u]) = f_n([u'])$, and let $(M,h)$,  $(M,h')$ be 
$\epsilon_0 \epsilon$--hermitian spaces over $D$
such that $[u], [u'] \in H^{\epsilon}(E,\sigma) $ correspond to the isomorphism classes
of $h$, respectively $h'$,  by the bijection of lemma  3.8.1. The
 elements $[u_n], [u_n'] \in H^{\epsilon}(E_n,\sigma_n)$
correspond to the isomorphism classes of $n \bullet h$, respectively $n \bullet h'$. 
We have $[u_n] = f([u]) = f([u']) = [u'_n]$, hence $n \bullet h \simeq n \bullet h'$. Since $k$ is real closed, corollary 3.6.2, applied to $D$ instead of $R$, shows that $h \simeq h'$. Therefore $[u] = [u']$, and this proves the
injectivity of $f_n$ in the case where $E$ is a simple algebra.

\smallskip
\n
{\it The case of a semi--simple algebra.}

Suppose now that $E$ is semi--simple. Then 

\quad $$E \simeq E_1 \times \dots \times E_r \times
A \times A^{{\rm op}}, $$

\n where $E_1, \dots, E_r$ are simple algebras which are stable under the involution
$\sigma$, and where the restriction of $\sigma$ to $A \times A^{{\rm op}}$ exchanges the two
factors. Applying lemma 3.7.1 and lemma 3.7.2 we are reduced to the case where
$E$ is a simple algebra, and we already know that the result is true in this case.

\smallskip
\n
{\it General case.}  

 Let $\overline E = E/{\rm rad}(E)$. Then $\overline E$ is semi--simple,
and $\sigma$ induces a $k$--linear involution $\overline \sigma : \overline E \to \overline E$.
Set $\overline E_n = E_n/{\rm rad}(E_n)$. Then $\overline E_n = {\bf M}_n(\overline E)$, 
and $\sigma_n$ induces $\overline \sigma_n : \overline E_n \to  \overline E_n$. We have
the following commutative diagram
$$\matrix {H^{\epsilon}(E,\sigma) & {\buildrel {f_n} \over \longrightarrow}  & H^{\epsilon}(E_n,\sigma) \cr 
\downarrow & {} & \downarrow \cr 
H^{ \epsilon}(\overline E, \overline \sigma) 
&{\buildrel {\overline f_n}  \over \longrightarrow}  & H^{ \epsilon}(\overline E_n,\overline \sigma), \cr }$$ 
where
the vertical maps are induced by the projection $E \to \overline E$. By lemma 3.7.3, these
maps are bijective. As $\overline E$ is semi--simple,  the map $\overline f_n$ is injective, 
hence $f_n$ is  also injective. This concludes the proof. 

\medskip

{\bf 3.11. Proofs of theorem 3.5.3 and theorem 3.1.1}

\smallskip

As we have already seen at the end of \S3.5, it is enough to prove theorem 3.5.3. This means proving that, if $k$ is real closed, and if $n > 0$, then
every $\epsilon$-hermitian space $(M,h)$ over $R$ has property ${\rm Div}_n$. Let $(M',h')$ be an $\epsilon$-hermitian space over $R$ such that $n\bullet (M,h) \simeq n\bullet (M',h')$. As noticed in \S3.1,  we may assume that $M' = M$.
Let $E = {\rm End}(M)$ and let $\tau$ be the involution of $E$ induced
by $h$, cf. \S3.8.  As explained in \S3.8 (applied to $h_0 = h$), the hermitian forms $h$ and $h'$ define $\epsilon$-hermitian elements $u = h_0^{-1}h = 1$ and $u' = h_0^{-1}h'$ of $(E,\tau)$; let $[u]$ and $[u']$ be their classes
in $H^\epsilon(E,\tau)$. Since $n\bullet (M,h) \simeq n\bullet (M,h')$, we have $[u_n] = [u'_n]$, where  $u_n$ and $u'_n$ are defined as in proposition 3.9.1 and theorem 3.10.1. By theorem 3.10.1 this implies $[u]=[u']$, hence  $(M,h) \simeq (M,h')$ by lemma 3.8.1. This concludes the proof.

\smallskip

\n {\it Question.} In several of the proofs above, we have made an essential use of the existence of a commutative subfield $k$ of the center of $R$ having  the 
following two properties:
 
 a) The involution of  $R$  is $k$-linear.
 
 b) ${\rm dim}_k R < \infty.$
 
\n We do not know what happens when there is no such $k$. For instance, suppose that
 $R$ is a skew field with involution which is of infinite dimension over its center, and that  $n$ 
 is odd; does ${\rm Div}_n$ hold for the $R$-hermitian spaces of finite dimension over $R$ ? It does\footnote{ Use the fact that, over $\F_p,  p \neq 2$, a quadratic form $<\!1,...,1\!>$ of odd rank is isomorphic to the direct sum of an hyperbolic form and 
 either $<\!1\!>$ or $<\!-1\!>$.} when the characteristic is $\neq 0$. But what about the characteristic 0 case ?
\medskip

\centerline
  {\bf \S 4. Induction and restriction for $G$-quadratic forms}

\medskip

   Let $k$ be a field of characteristic $\not= 2$. Let $G$ be a finite group and let $S$ be a 2-Sylow subgroup of $G$. Our aim is to prove th.1.1.2.
   We recall its statement:

   \medskip
\n
{\bf Theorem 1.1.2} - {\it Let $(V_1,q_1)$ and $(V_2,q_2)$ be two $S$-quadratic spaces over $k$. Suppose that}:

\smallskip

 $${\rm Res}^G_S \ {\rm Ind}^G_S(V_1,q_1) \simeq_S{\rm Res}^G_S \ {\rm Ind}^G_S(V_2,q_2).$$

\smallskip

\n {\it  Then}:
$${\rm Ind}^G_S(V_1,q_1) \simeq_G  {\rm Ind}^G_S(V_2,q_2).$$

\smallskip

\n
The proof is based on:

\smallskip
$\bullet$ the Witt-type cancellation theorem, cf. \S 2;

\smallskip
$\bullet$ the ``no odd torsion" principle, cf. \S 3;

\smallskip
$\bullet$ the elementary properties of the Burnside rings,
cf. \S4.2 below.
    
    \medskip \medskip
    
    {\bf 4.1. Algebraic properties of the projection formula}
    
\medskip
  Let $A$ and $B$ be two commutative rings, and let $r: B\to A$ and $i:A\to B$
  be respectively a ring homomorphism and an additive homomorphism
  such that:
  
  \smallskip
  (4.1.1) \hskip .5 cm  $i(r(y)x) = y.i(x)$  for every  $x\in A$ and $y\in B$.
  
 \smallskip
\n   [This formula means that $i: A \to B$ is $B$-linear if $A$ is viewed as a
  $B$-module via $r$. In $K$-theory, such a formula is often called a ``projection formula'', and  $i$ and $r$ are denoted by symbols such as $f_*$ and $f^*$. In \S 4.3 below,  $i$ will be ${\rm Ind}^G_S$ and
  $r$ will be ${\rm Res}^G_S$.]
  
  \medskip

\n We put

\medskip

(4.1.2) \hskip .4 cm $Q = i(1_A) \in B$ \hskip .4 cm  and \hskip .4 cm  $q = r(Q) \in A,$

\smallskip

\n where $1_A$ denotes the unit element of the ring $A$.

\medskip

\n Formula (4.1.1), applied to $x=1_A$, gives

\medskip

(4.1.3) \hskip .5 cm  $ i(r(y)) = Q.y$ \  \ for all $y\in B$.

\medskip
\n Let us define $R : A \to A$ \ by:

\smallskip

(4.1.4) \quad $R = r \circ i.$

\smallskip

\n It is an additive endomorphism of $A$; we want to compare its kernel to that of $i$. To do so, let $(F,n)$ be a pair where  $F \in {\bf Z}[t]$
is a one-variable polynomial with coefficients in  $\bf Z$ and  $n$  is an integer; assume that:
\smallskip

(4.1.5) \hskip .5 cm $q.F(q) = n.1_A$.

\medskip

\n {\bf Proposition 4.1.1 } - {\it If $A,B,i,r,F,n$ are as above, we have} :

\smallskip

(4.1.6)  \hskip 1cm   $n.i(a) = i(F(q).R(a))$ \ {\it for every} \ $a \in A$.

\medskip

\n  {\it Proof.} Let $q' = F(q)$ and let $Q' =F(Q)$. We have:

\smallskip

\n $i(q'.R(a))=Q'.i(R(a))$ by (4.1.1) applied to $ x=R(a)$ and \ $y=Q'$

\smallskip
\ \ \ \qquad $= Q'Q.i(a)$ by (4.1.3) applied to $y=i(a)$

\smallskip

\ \ \ \qquad $=i(r(Q'Q)a)$ by (4.1.1) applied to $x=a$ and \ $y=Q'Q$

\smallskip

\ \ \ \qquad $=i(n.a) $ since $r(Q'Q)=q'q=n.1_A$

\smallskip

\ \ \ \qquad $= n.i(a)$.

\medskip\medskip
 
 \n  {\bf Corollary 4.1.2} - {\it Assume that $B$ is $n$-torsion-free. Then ${\rm Ker} (i) =  {\rm Ker} (R).$}

\smallskip \n  {\it Proof.} The inclusion ${\rm Ker}(i) \subset {\rm Ker}(R)$ is obvious, and the opposite inclusion follows from (4.1.6), which shows that
$n.i(a) = 0$ if $R(a)=0$.

\smallskip
\n  {\it Remark.} Corollary 4.1.2 is all we need for the proof of theorem 1.1.2.
However, the following result is worth mentioning:

\smallskip

\n  {\bf Proposition 4.1.3} - (1) {\it The kernel and the cokernel of the natural map

\smallskip
\hskip 4cm   ${\rm Ker}(R) \oplus {\rm Im}(R)  \  \to  \  A $

\n  are killed by $n$.}

\smallskip
\n  (2) {\it The same is true for the map \ ${\rm Ker}(r) \oplus {\rm Im}(i) \to B.$}

\smallskip
\n [Recall that an abelian group  $X$ is {\it killed by $n$} if $n.x=0$ for every $x\in X.$]

\medskip

\n  {\it Proof of} (1). We have

\smallskip
(4.1.7) $ R^2(a) = q.R(a)$ \ for every $a \in A$.

\smallskip
\n Indeed, by (4.1.3) we have $i(r(i(a)))=Q.i(a)$, hence

\smallskip
$R^2(a)  = r(i(r(i(a)))) = r(Q.i(a)) = r(Q).r(i(a)) = q.R(a).$

\smallskip
Let now $x$ be an element of ${\rm Ker}(R) \cap {\rm Im}(R)$, and let us write $x$
as $R(a)$, with $a \in A$. By (4.1.7), we have $q.R(a) = 0$, i.e. $q.x=0$,
hence $n.x$ = 0 by (4.1.5). This shows that the kernel of ${ \rm Ker}(R) \oplus {\rm Im}(R)  \  \to  \  A $ is killed by $n$. In order to prove the same result for the cokernel, we need to show that, for every $x\in A$, there exists $y \in {\rm Im}  (R)$
and $z\in {\rm Ker}(R)$ with $nx = y+z$. We choose $y = F(q).R(x)= r(F(Q).i(x))= R(F(q).x)$, and we define $z$ as $z = n.x-y$. It remains to check that $R(z)=0$, i.e. that $n.R(x)=R(y)$; this follows from $R(y) = R^2(F(q).x) = q.R(F(q).x)=
q.F(q).R(x) = n.R(x)$.

\medskip

\n  {\it Proof of } (2). If $b\in B$ belongs to ${\rm Ker}(r) \cap {\rm Im}(i)$, we may write $b$ as $i(a)$ with $a \in A$, and we have $r(i(a))=0$, i.e. $R(a)=0$, which implies $n.i(a) = 0$ by (4.1.6), i.e.  $n.b = 0$.
 \smallskip 
   If $b$ is any element of $B$, put $y=Q.F(Q).b$ and $z=n.b-y$. We have
   $r(y)=q.F(q).r(b) = n.r(b)$, hence $r(z)=0$. Since $y$ belongs to
   ${\rm Im}(i)$ by (4.1.3), the relation $n.b = y + z$ shows that $n.b$ belongs to
${\rm  Ker}(r) \oplus {\rm Im}(i) .$

\medskip

\n  {\it Remark.} The above results are only interesting when the integer $n$ can be chosen such that
the rings $A$ and $B$ have no $n$-torsion. In the case we need, it will be enough that $n$ be {\it odd}; this is what the Burnside ring method is going to give us.

\medskip

 {\bf 4.2. The Burnside ring}
    
    \medskip
    
    \n  {\it Definitions and Notation}
    
    \smallskip
    
    If $G$ is a finite group, we denote by ${\rm Burn}(G)$ its Burnside ring, i.e.
    the Grothendieck group of the category of finite $G$-sets, cf. [B 91, \S5.4]. Every finite $G$-set $X$ defines an element $b_X$ of ${\rm Burn}(G)$,
    and we have:
    
    \smallskip
     
      (4.2.1)  \hskip 1cm   $b_{X \sqcup Y} = b_X + b_Y$

      \n  and
      
      (4.2.2) \hskip 1cm  $b_{X \times Y} = b_X.b_Y.$
      
      \smallskip
      
      \n  Let $\cal H$ be a set of representatives of the subgroups of $G$ up to $G$-conjugation. The family $\{b_{G/H}\}_{H\in \cal H}$ is a $\bf Z$-basis of ${\rm Burn}(G)$; hence ${\rm Burn}(G)$ is a free $\bf Z$-module of rank $h=|\cal H|$.
      If $H$ is a subgroup of $G$, there is a unique ring homomorphism
      
       \smallskip
       
      (4.2.3) \hskip 1cm  $f_H : {\rm Burn}(G) \ \to \ {\bf Z}$
      
       \smallskip
      
      \n  such that $f_H(b_X) = \ $number of points of $X$ which are fixed by $H$.
      
      The family $(f_H)_{H \in \cal H}$ gives a ring homomorphism
      
       \smallskip

      (4.2.4) \hskip 1cm  $ f: {\rm Burn} (G) \ \to \ {\bf Z}^{\cal H}$,
      
       \smallskip
       
       \n  where ${\bf Z}^{\cal H}$ denotes the product of $h$ copies of $\bf Z$ indexed
       by $\cal H$. This homomorphism is injective, and its cokernel is finite; 
       we may thus view ${\bf Z}^{\cal H}$ as the {\it normalization} of the ring ${\rm Burn}(G)$. 
       
   \n [As shown by A. Dress, this leads to an explicit description
       of ${\rm Spec \ Burn}(G)$, as a union of $h$ copies of ${\rm Spec} \ {\bf Z}$ with suitable glueing; the precise statement is given in  [B 91, th.5.4.6] where however $f_K(x)$ is misprinted as $f_H(x)$. A remarkable corollary is that ${\rm Spec \ Burn}(G)$ is connected if and only if $G$ is solvable, cf. [B 91, cor.5.4.8].]
      
      \medskip
      
      \n  {\it Characteristic polynomial and norm}
      
      \smallskip
      
      Let  $x$  be an element of ${\rm Burn}(G)$. The multiplication by $x$ is an endomorphism of the additive group ${\rm Burn}(G) \simeq {\bf Z}^h$. Let  $P_x(t) \in {\bf Z}[t]$ be its characteristic polynomial; it is a monic polynomial of degree $h$.
      By the Cayley-Hamilton formula, we have
      
      \smallskip
      
      (4.2.5) \hskip 1cm   $ P_x(x) = 0$  \ in  \ ${\rm Burn}(G)$.
      
      \smallskip
      
      The embedding (4.2.4) shows that $P_x(t)$ is also the characteristic polynomial of $x$ acting by multiplication on ${\bf Z}^{\cal H}$, hence :
      
      \smallskip
      
      (4.2.6) \hskip 1cm  $P_x(t) = \prod_{H\in {\cal H}} (t - f_H(x)).$
      
      \smallskip
      
      \n  The {\it norm}  $N(x)$ of $x$ is 
      
      \smallskip
      
      (4.2.7) \hskip 1cm   $N(x) = \prod_{H\in {\cal H}} f_H(x)$.
      
      \smallskip
      
      \n  It is the constant term of \ $(-1)^hP_x(t)$.
      
      \smallskip
      

 \medskip
      
      \n  {\bf Proposition 4.2.1} - {\it Let $p$ be a prime number, and suppose
      that $G$ is a $p$-group. Let  $X$ be a finite $G$-set and let $x = b_X$ be the corresponding element of ${\rm Burn} (G)$. Assume that $|X| = {\rm Card} X$ is not divisible by $p$. Then} :
      
 \smallskip     
      (i) $N(x)$ {\it is not divisible by} $p$.
      
\smallskip      
      (ii) {\it There exists  $F_x \in {\bf Z}[t]$ such that $x.F_x(x) = N(x)$ in}     ${\rm Burn}(G)$.  
      
      \medskip
      
      \n  {\it Proof}. Let $H$ be a subgroup of $G$ and let $X^H$ be the set of
      all elements of $X$ which are fixed by $H$. Since $H$ is a $p$-group,
      we have $|X^H| \equiv |X| \pmod {p}$, hence $f_H(x)= |X^H|$ is not divisible by $p$; by (4.2.7), the same is true for $N(x)$. This proves (i).
   \smallskip   
      Since $N(x)$ is the constant term of $(-1)^hP_x(t)$, we may define
      a polynomial $F_x(t)$ by the formula 
      
      \smallskip
      
       \hskip 1cm  $t.F_x(t) = N(x) - (-1)^hP_x(t)$,  
      
      \smallskip
      
      \n  and we have  $x.F_x(x) = N(x) - (-1)^hP_x(x) =  N(x)$, cf.  (4.2.5).  This proves (ii).

\medskip

   {\bf 4.3. Proof of theorem 1.1.2}
    
    \smallskip

   Let  $\Lambda(G)$ denote the Grothendieck-Witt group 
 of the $G$-quadratic spaces over $k$, i.e. the group denoted by ${\rm {\rm WGr}}^1(k[G])$ in \S2.6. The tensor product operation makes $\Lambda(G)$ a {\it commutative ring} [this reflects the fact that $k[G]$ is a {\it bigebra} and not merely an algebra]. 
 
   The ring $\Lambda(S)$ is defined in a similar way. The induction and
   restriction functors give rise to maps
   
   \smallskip
   
     ${\rm Ind}^G_S: \Lambda(S) \to  \Lambda(G)$        \ and \ 
     ${\rm Res}^G_S :  \Lambda(G)   \to \Lambda(S)$.
     
     \smallskip
     
  \n   The rings $A = \Lambda(S)$ and $B = \Lambda(G)$, together with the maps
     $i = {\rm Ind}^G_S$ and $r= {\rm Res}^G_S,$ have the properties described at the beginning of \S4.1; in particular, formula (4.1.1) is valid, as a simple computation shows. We may then apply the constructions of \S4.1, and define:
     
     \smallskip
     
       \quad $Q = i(1_A) \in B$ \quad and \quad  $q = r(Q) \in A.$

\medskip

\n  {\bf Lemma 4.3.1} - {\it There exists 
$F \in {\bf Z}[t]$ and $n \in {\bf Z}$, with $n$ odd, such that $q.F(q) = n.1_A$, as in $(4.1.5)$.}

\smallskip

\n  {\it Proof}. Let ${\rm Burn}(G)$ be the Burnside ring of the group  $G$,
 cf \S4.2. There is a unique ring homomorphism

\smallskip

 \quad   $\varepsilon_G : {\rm Burn}(G) \ \to \ B =  \Lambda(G)$,

\smallskip

\n  having the following property: if $X$ is any finite $G$-set, $\varepsilon_G(b_X)$ is equal to $[q_X]$ , where  $q_X$ is the unique $G$-quadratic form with an orthonormal basis isomorphic to $X$ as a $G$-set.

\smallskip

  Let now  $X=G/S$, with the natural action of $G$. The corresponding $G$-quadratic form is $
  {\rm Ind}^G_S({\bf 1})$ , where $\bf 1$ means
  the unit quadratic form of rank 1, with trivial action of $S$; its class in $B$ is  $Q = i(1_A)$.
  
  \smallskip
  We may apply the same construction to $S$; we have
  
  \smallskip
  
    \quad $\varepsilon_S :  {\rm Burn} (S) \ \to \ A =  \Lambda(S)$,
  
  \smallskip
  
\n  and $\varepsilon_S(b_X)= r(Q) = q$. By proposition 4.2.1, applied to $S$, to $ X=G/S,$
  and to $p=2$, there exist $F \in {\bf Z}[t]$ and $n \in {\bf Z}$, with $n$ odd, such that
  $b_X.F(b_X) = n$ in ${\rm Burn}(S).$ Applying $\varepsilon_S$ to this equality gives  $q.F(q) = n.1_A$ in $A$, as wanted.
  
  \medskip
  
  \n  {\it Remark}. Using the Burnside ring in questions involving
  $G$-functors is a very efficient technique, which was introduced by
  A. Dress in the early 1970s. See e.g. [B 91, \S 5.6].
  
  \smallskip

\medskip

\n  {\bf Theorem 4.3.2} - {\it The maps

\smallskip
\qquad $ {\rm Res}^G_S \circ {\rm Ind}^G_S :  \Lambda(S)   \to \Lambda(S)$  \ \ and \ \
 ${\rm Ind}^G_S: \Lambda(S) \to  \Lambda(G)$ 
 
 \smallskip
\n have the same kernel.}

 \smallskip
 
 \n {\it Proof.} Let $n$ be as in lemma 4.3.1.  Since $n$ is odd,  the ring $B =  \Lambda(G)$ is $n$-torsion-free (see theorem 3.1.1). By corollary 4.1.2, we have ${\rm Ker}(i) = {\rm Ker}(r \circ i)$.
     
\medskip
  
  \n  {\it End of the proof of theorem 1.1.2.} 
  
   Let $v_1$ and $v_2$ be the classes of $(V_1,q_1)$ and $(V_2,q_2)$ in  $\Lambda(S)$. By assumption, $v_1-v_2$ belongs to the kernel of $ {\rm Res}^G_S \circ {\rm Ind}^G_S$; by theorem 4.3.2, it belongs to the kernel of  $ {\rm Ind}^G_S$.
  By  formula (2.6.1), this means that
 ${\rm Ind}^G_S(V_1,q_1)$ is $G$--isomorphic to $ {\rm Ind}^G_S(V_2,q_2).$
  \bigskip
  
\centerline
{\bf References}

\medskip

\medskip
[A6] N. Bourbaki, Alg\`ebre, Chap.6, {\it Groupes et corps ordonn\'es}, Masson,
Paris, 1981; new printing, Springer--Verlag, 2006; English translation, Springer-Verlag, 1998. 

\smallskip

[A8]   \ -------- \ , Alg\`ebre, Chap.8,
{\it Modules et anneaux semi--simples}, new revised edition, Springer--Verlag, 2012. 

\smallskip
[A9]   \ -------- \ , Alg\`ebre, Chap.9, {\it Formes sesquilin\'eaires et formes
quadratiques}, Hermann, Paris, 1959; new printing, Springer-Verlag, 2006.

\smallskip
[BFL 90] E. Bayer--Fluckiger \& H.W. Lenstra, Jr.,  {\it Forms in odd degree
extensions and self-dual normal bases},
Amer. J. Math. {\bf112} (1990), 359-373.

\smallskip
[BFP 11] E. Bayer--Fluckiger \& R. Parimala, 
{\it Galois Algebras, Hasse Principle and Induction--Restriction Methods},
Documenta Math. {\bf 16} (2011), 677--707; {\it Errata}, to appear. 

\smallskip

[BFS 94] E. Bayer--Fluckiger \& J--P. Serre, {\it Torsions quadratiques
et bases normales autoduales}, 
Amer. J. Math. {\bf116}  (1994),
1--64.

\smallskip
 [B 91] D.J. Benson, {\it Representations and Cohomology I}, Cambridge studies in advanced mathematics  {\bf 30}, Cambridge University Press, Cambridge, 1991.
 
\smallskip
[K 91] M. Knus, {\it Quadratic and Hermitian Forms over Rings}, Grundlehren math. Wiss. {\bf294}, Springer-Verlag, 1991.

\smallskip
[KMRT 98] M. Knus, A. Merkurjev, M. Rost \& J--P. Tignol, {\it The Book of Involutions}, AMS Colloquium Publications {\bf 44}, 1998. 



\smallskip
[M 86] J. Morales, {\it Integral bilinear forms with a group action},  J. Algebra {\bf98} (1986), 470-484. 

\smallskip
[QSS 79] H--G. Quebbemann, W. Scharlau \& M. Schulte,  {\it Quadratic and hermitian
forms in additive and abelian categories}, J. Algebra  {\bf59} (1979), 264--289.

\smallskip
[QSSS 76] H--G. Quebbemann, R. Scharlau, W. Scharlau \& M. Schulte, {\it Quadratische
Formen in additiven Kategorien},   Bull. Soc. Math. France, M\'emoire {\bf48} (1976), 93--101.

\smallskip
[R 11] C.R. Riehm, {\it Introduction to Orthogonal, Symplectic and Unitary Representations of Finite Groups}, Fields Institute Monographies 28, AMS, 2011.


\smallskip
[S 85] W. Scharlau, {\it Quadratic and Hermitian Forms}, Grundlehren math. Wiss. {\bf 270},
Springer--Verlag, 1985.

\medskip
\medskip

E. Bayer--Fluckiger

\'Ecole Polytechnique F\'ed\'erale de Lausanne

EPFL/FSB/MATHGEOM/CSAG

Station 8

1015 Lausanne, Switzerland

\smallskip
eva.bayer@epfl.ch

\medskip
\medskip
R. Parimala

Department of Mathematics $ \&$ Computer Science

Emory University

Atlanta, GA 30322, USA

\smallskip
parimala@mathcs.emory.edu

\medskip
\medskip
J--P. Serre

Coll\`ege de France

3 rue d'Ulm

75005 Paris, France

\smallskip
serre@noos.fr

\end{document}